\documentclass[12pt]{article}

\usepackage[T2A]{fontenc}
\usepackage[utf8]{inputenc}

\usepackage{amsfonts}
\usepackage{amssymb}
\usepackage{amsmath}
\usepackage{amsthm}

\def \le {\leqslant}
\def \ge {\geqslant}

\begin{document}
\begin{large}
\centerline{A note on well distributed sequences}
\end{large}
 \vskip+0.3cm

\centerline{Nikolay Moshchevitin}
 \vskip+1cm

{\bf 1. Well distributed sequences.}
Let 
\begin{equation}\label{1}
\pmb{\xi}_k=  (\xi_{1,k},...,\xi_{n.k})   \,\,\,\, k=1,2,3,...
\end{equation}
be an infinite  sequence of points in the unit cube $[0,1)^n$. 
For $\pmb{\eta} = (\eta_1,...,\eta_n) \in [0,1)^n$ we consider boxes of the form
\begin{equation}\label{box}
B_{C,q}(\pmb{\eta})=
\left[\eta_1, \eta_1 + \frac{C}{q^{1/n}}\right)\times \cdots\times   \left[\eta_n, \eta_n + \frac{C}{q^{1/n}}\right) \subset [0,1)^n.
\end{equation}
We define  sequence (\ref{1}) to be {\it well distributed} if there exists a positive constant $C$  and infinitely many positive integers $ q$ such that 
for any box of the form (\ref{box})  there exists  a positive integer $ k \le q$ with 
\begin{equation}\label{D}
 \pmb{\xi}_k \in  B_{C,q}(\pmb{\eta}).
\end{equation} 
Chebyshev showed that for any irrational $\alpha$  there exists infinitely many positive integers $q$ such that for any $\eta \in [0,1)$ there exists
positive integer $k \le q$ such that 
$$
||k \alpha -\beta ||\le \frac{3}{q}.
$$  
This means that the sequence
\begin{equation}\label{dd}
\{ k\alpha \} ,\,\,\,  k =1,2,3,...
\end{equation}
 is well distributed.
Generalising a result by Chebyshev, Khintchine \cite{H} proved that the Kronecker sequence
\begin{equation}\label{2}
\pmb{\xi}_k = (\{ k\alpha_1\},...,\{ k\alpha_n\}) \,\,\,\, \alpha_1,...,\alpha_n \in \mathbb{R}
\end{equation}
is well distributed if and only if vector $\pmb{\alpha} = (\alpha_1,...,\alpha_n)\in \mathbb{R}^n$ is non-singular.
For the sake of completeness in the next section we formulate this result as it is given in Cassels' book \cite{C}.
Certain discussion about different formulations of this result one can find in \cite{M}. Here we would like to add that Khintchine's result is valid not only for the Kronecker sequence  but for general systems of $n$ linear forms in $m$ variables.

 \vskip+0.3cm

{\bf 2. Khintchine's criterium.}
Consider irrationality measure function
$$
\psi_{\pmb{\alpha}} (t) = \min_{q\in \mathbb{Z}_+: \, q\le t }\, \max_{1\le j\le n} || q\alpha_j||,\,\,\,\,\,\,\, ||\cdot || = \min_{a\in \mathbb{Z}}|\cdot - a|.
$$
By Dirichlet theorem 
$$
\psi_{\pmb{\alpha}} (t) \le \frac{1}{t^{1/n}}\,\,\,\,\,\forall t\ge 1.
$$
Vector $\pmb{\alpha} \in \mathbb{R}^n$ is called {\it non-singular} if
$$
\limsup_{t\to\infty}  t^{1/n}
\psi_{\pmb{\alpha}} (t) >0.
$$
Theorem XIII from \cite{C}  (to be more precise  formula (4) from the proof of Theorem XII which is based on Theorem VI) together with transference theorem for singular vectors (Theorem XII) states that $\pmb{\alpha}$ is non-singular if and only if
$$
\exists \, \delta >0\,\,\,
\text{and}\,\,\,
\exists \,\, \text{infinitely many}\,\, q\in \mathbb{Z}_+\,\,\text{such that}\,\,
\forall \pmb{\eta} = (\eta_1,...,\eta_n)\in \mathbb{R}^n
$$
$$
\exists k  \in \mathbb{Z}_+:\, k \le q, \,\, \max_{1\le j \le n} ||k\alpha_j- \eta_j||\le \frac{\delta}{q^{1/n}}. 
$$
The last statement means that  for $C = 2 \delta$ in any box of the form 
 (\ref{box})
there exists a point of the sequence (\ref{1}) with $ k \le q$.

 \vskip+0.3cm

{\bf 3. General statement.}
In this short communication we prove the following  result.

 \vskip+0.3cm
 {\bf Theorem 1.} \, {\it
Let  sequence  (\ref{1}) be well distributed. Then for almost all $\pmb{\eta} = (\eta_1,...,\eta_n) \in [0,1)^n$
one has
$$
\liminf_{k \to \infty} k^{1/n}
\max_{1\le j \le n} || \xi_{j,k} - \eta _j|| =0
. 
$$
}

 \vskip+0.3cm

{\bf Remark.} For the special case of  sequence (\ref{dd}) this result coincide with a statement from  \cite{D}, which was proved there  by means of Ostrowski  numerical systems.
Multidimensional results which generalise one-dimensional theorem from \cite{D} and deal with systems of $n$ linear forms in $m$ variables were  obtained in \cite{U}. The results from
\cite{U} deal with a setting in terms of Dynamical systems on space of lattices. Our definition of well distribution is close to the  definition fo
local ubiquity from \cite{B1}.   Another proof of a statement similar to our Theorem 1 can be deduced from Lemma 1 from \cite{B2}.

 \vskip+0.3cm
Proof of Theorem 1.
We take positive reals $\psi_\nu< 1, \lim_{\nu \to \infty} \psi_\nu  =0$ such that the series $\sum_{\nu=1}^\infty \psi_\nu^n $ diverges and small $\varepsilon >0$. Below $\mu (\cdot )$ stands for the Lebesgue measure.

 By the conditions of Theorem 1  there exists  $C>0$  and arbitrary large integers $q$ such that 
 in any box of the form  (\ref{box}) there exists a point of the sequence (\ref{2})
with $ k \le q$. 
We cover the cube $[0,1)^n$ by  $ W = \left\lceil \frac{q^{1/n}}{C}\right\rceil^n$  boxes 
$B_{C,q}(\pmb{\eta}_{i_1,...,i_n})$
of the form (\ref{box})
with 
$$\pmb{\eta}_{i_1,...,i_n} = \left(\frac{i_1C}{q^{1/n}},..., \frac{i_nC}{q^{1/n}}\right), \,\,\,\, 0\le i_1,...,i_n <  \left\lceil \frac{q^{1/n}}{C}\right\rceil^n,\,\,\, i_1,...,i_n \in \mathbb{Z}
$$
 which do not intersect  by inner points. 
At least of $W' = \left( \left\lceil \frac{q^{1/n}}{C}\right\rceil -1\right)^n$ of these boxes lies inside the cube $[0,1)^n$.
By well approximability, in each of these last boxes there is a point of the form (\ref{1}) with
$ 1\le k \le q.
$
For positive $ \psi <1$  we take $q'$ boxes of the form
 \begin{equation}\label{l}
I_l=  \left[\xi_{1,k_l}- \frac{C\psi}{q^{1/n}},\xi_{1,k_l}+\frac{C\psi}{q^{1/n}}\right]\times\cdots\times
 \left[\xi_{n,k_l}-\frac{C\psi}{q^{1/n}},\xi_{n,k_l}+\frac{C\psi}{q^{1/n}}\right],\,\,\,\,1\le l \le q'
 \end{equation}
 with the  centres at certain points $\pmb{\xi}_{k_l}  , \, 1\le k_l \le q$ 
 which belong to boxes 
 \begin{equation}\label{Cq}
 B_{C,q}(\pmb{\eta}_{i_1,...,i_n})\,\,\,\,\,\text{with} \,\,\,\,\,i_1\equiv ...\equiv i_n \equiv 1 \pmod{3}. 
 \end{equation}
 In each of such boxes $B_{C,q}(\pmb{\eta}_{i_1,...,i_n})$ we take just one point  $\pmb{\xi}_{k_l} $.
Then
$$
I_l\cap I_{l'}=\varnothing,\,\,\,l\neq l'.    
 $$
and
$
 q'\asymp_C q.
$

Now from these integers $q, q'$  we construct  sequences $q_\nu$  and $q_\nu'$  by an inductive  procedure.
As $q_1$ we take an integer $q$ from the definition of the well distribution. It  correspond to a certain $q_1'$ and  can be chosen to be arbitrary large.

 Now we suppose that the numbers $ q_1,  ..., q_\nu$ satisfying the definition of the well approximability are defined as well as the  numbers 
 $ q_1',  ..., q_\nu'$ satisfying 
  \begin{equation}\label{c1}
 q_\nu'\asymp_C q_\nu.
 \end{equation}
 and the boxes
  \begin{equation}\label{l1}
I_l(\nu)  =  \left[\xi_{1,k_l^\nu}- \frac{C\psi_\nu}{q_\nu^{1/n}},\xi_{1,k_l^\nu}+\frac{C\psi_\nu}{q_\nu^{1/n}}\right]\times\cdots\times
 \left[\xi_{n,k_l^\nu}- \frac{C\psi_\nu}{q_\nu^{1/n}},\xi_{n,k_l^\nu}+\frac{C\psi_\nu}{q_\nu^{1/n}}\right],\,\,\,\,1\le l \le q_\nu',
 \end{equation}
 where $1\le k_l^\nu \le q_\nu$ and 
 \begin{equation}\label{cc1}
I_l(\nu)\cap I_{l'}(\nu)=\varnothing,\,\,\,l\neq l'.    
 \end{equation}
 We can chose $q=q_{\nu+1}$ from the definition of well aproximabolity   large enough to satisfy
 \begin{equation}\label{sharp}
 \sharp \{ r:\, 1\le r\le q_{\nu+1}',\,\,\,  I_r(\nu+1) \cap I_l(\lambda) \neq \varnothing \}  \le (1+\varepsilon)
{\mu (I_l(\lambda))}{q_{\nu+1}'}   \end{equation}
 for every box $I_l(\lambda), 1\le l \le q_\lambda', 1\le\lambda \le \nu$  which appeared during the  previous steps $\lambda \le \nu$.
 This is possible by the definition of well approximability because the number of boxes of the form (\ref{Cq}) with $q=q_{\nu+1}$ which intersect  given box
  $I_l(\lambda) $
 is
 $
 \sim
 {\mu (I_l(\lambda))}{q_{\nu+1}'} $ when $ q=q_{\nu+1} \to \infty$.

  Now we define the union
 $$
 E_\nu = \bigcup_{l=1}^{q_\nu'} I_l(\nu).
 $$
 First of all by (\ref{c1}) and (\ref{cc1}) we see that 
 $$
 \mu (E_\nu)=  \sum_{l} \mu (   I_l(\nu))\gg_C \psi_\nu^n
 $$
 and so
  \begin{equation}\label{ccc1}
  \sum_{\nu=1}^\infty
 \mu (E_\nu) = \infty.
 \end{equation}
 Then for $ \lambda < \nu$ we have
   \begin{equation}\label{ccc2}
 \mu (E_\lambda \cap E_\nu) \le (1+\varepsilon)  \mu (E_\lambda)  \mu (E_\nu)
.
 \end{equation}
 A famous lemma from metric theory of Diophantine approximation 
 (see \cite{S}, Lemma 5 \S 3 Ch. 1) shows that the measure of the set 
 $$
  E = \{ \pmb{\eta}:\,\, \exists \,\,\text{infinitely many}\,\, \nu\,\,
  \text{such that }\,\, \pmb{\eta} \in E_\nu\}
  $$
  is
  $$
  \ge
  \limsup_{t\to \infty}
  \frac{\left(\sum_{\nu=1}^t   \mu (E_\nu) \right)^2}{  \sum_{\lambda,\nu =1}^t \mu (E_\lambda \cap E_\nu)} \ge \frac{1}{1+\varepsilon},
  $$
  by (\ref{ccc1}) and (\ref{ccc2}). As $\varepsilon $ is arbitrary, Theorem 1 in proven.
 
  \vskip+0.3cm
 {\bf Acknowledgements.}
 The research was funded by the Russian Science Foundation (project No. 22-41-05001).

author: Nikolay Moshchevitin

Big Data and Information Retrieval School

Faculty of Computer Science

National Research University Higher School of Economics

and

Israel Institute of Technology (Technion)

emails: moshcheviitin@gmail.com, moshchevitin@technion.ac.il

\end{document}